\documentclass[a4paper,11pt]{article}
\usepackage{amsfonts,amssymb,amsthm,cite,amsmath,amstext}
\usepackage{color}
\usepackage[colorlinks,linkcolor=black,citecolor=black]{hyperref}
\usepackage{enumerate}
\numberwithin{equation}{section}

\usepackage{amsmath,amssymb,amsthm}
\usepackage{amsmath,amssymb,amsthm,amscd}
\usepackage{color}
\usepackage[mathscr]{eucal}

\DeclareMathOperator{\vol}{vol}
\title{\textbf{A remark on the convergence of inverse $\sigma_k$-flow}}
\author{\textsc{Jian Xiao}}
\date{}
\begin{document}
\maketitle

\theoremstyle{definition}
\newtheorem*{pf}{Proof}
\newtheorem{theorem}{Theorem}[section]
\newtheorem{remark}{Remark}[section]
\newtheorem{problem}{Problem}[section]
\newtheorem{question}{Question}[section]
\newtheorem{conjecture}{Conjecture}[section]
\newtheorem{lemma}{Lemma}[section]
\newtheorem{corollary}{Corollary}[section]
\newtheorem{definition}{Definition}[section]
\newtheorem{proposition}{Proposition}[section]
\newtheorem{example}{Example}[section]
\newtheorem{observation}{Observation}[section]

\begin{abstract}
We study the positivity of related cohomology classes
concerning the convergence problem of inverse $\sigma_k$-flow in the conjecture proposed by Lejmi and Sz\'{e}kelyhidi.
\end{abstract}

\section{Introduction}
We aim to study the positivity of related cohomology classes in the following conjecture proposed by Lejmi and Sz\'{e}kelyhidi \cite{lejmisze15Jflow}. We generalize their conjecture by weakening the numerical condition on $X$ a little bit.

\begin{conjecture}
\label{conj LS15 sigma k flow}
(see \cite[Conjecture 18]{lejmisze15Jflow})
Let $X$ be a compact K\"ahler manifold of dimension $n$, and let $\omega,\alpha$ be two K\"ahler metrics over $X$ satisfying
\begin{align}
\label{eq over X}
\int_X \omega^n-\frac{n!}{k!(n-k)!}\omega^{n-k}\wedge \alpha^k \geq 0.
\end{align}
Then there exists a K\"ahler metric $\omega' \in \{\omega\}$
such that
\begin{align}
\label{eq n-1 positive}
{\omega'} ^{n-1}- \frac{(n-1)!}{k!(n-k-1)!}{\omega'} ^{n-k-1}\wedge \alpha^k >0
\end{align}
as a smooth $(n-1,n-1)$-form if and only if
\begin{align}
\label{eq numerical positive}
\int_V \omega^p - \frac{p!}{k!(p-k)!} \omega^{p-k}\wedge \alpha^k >0
\end{align}
for every irreducible subvariety of dimension $p$ with
$k\leq p \leq n-1$.
\end{conjecture}

For the previous works closely related to this conjecture, we refer the reader to \cite{Donaldson99Jflow_momentmap}, \cite{Chen00mabuchi_Jflow, Chen04Jflow}, \cite{songweinkove08Jflow} and \cite{maxinan11sigmaflow}.
And in this note we mainly concentrate on the case when $k=1$ and $k=n-1$. 

For $k=1$, \cite[Theorem 3]{collins2014convergence-Jflow} confirmed this conjecture for toric manifolds. Over a general compact K\"ahler manifold, it is not hard to see the implication (\ref{eq n-1 positive}) $\Rightarrow$ (\ref{eq numerical positive}) holds. In the reverse direction, we prove  $\{\omega-\alpha\}$ must be a K\"ahler class under the numerical conditions in Conjecture \ref{conj LS15 sigma k flow} for $k=1$; indeed, this is a necessary condition of (\ref{eq n-1 positive}) and \cite[Proposition 14]{lejmisze15Jflow} proved this over K\"ahler surfaces.

\begin{theorem}
\label{thm kaehler}
Let $X$ be a compact K\"ahler manifold of dimension $n$, and let $\omega,\alpha$ be two K\"ahler metrics over $X$ satisfying the numerical conditions in Conjecture \ref{conj LS15 sigma k flow} for $k=1$. Then $\{\omega-\alpha\}$ is a K\"ahler class.
\end{theorem}

For $k=n-1$, we have the following similar result.

\begin{theorem}
\label{thm n-1 flow}
Let $X$ be compact K\"ahler manifold of dimension $n$, and let $\omega, \alpha$ be two K\"ahler metrics over $X$ satisfying
the numerical conditions in Conjecture \ref{conj LS15 sigma k flow} for $k=n-1$. Then the class $\{\omega^{n-1}-\alpha^{n-1}\}$ lies in the closure of the Gauduchon cone, i.e. it has nonnegative intersection number with every pseudoeffective $(1,1)$-class.
\end{theorem}


\section*{Acknowledgements}
I would like to thank Philippe Eyssidieux and Mehdi Lejmi for introducing this problem to me, and thank Tristan C. Collins for useful comments.
This work is supported by the China Scholarship Council.

\section{Proof of the main results}
In this section, we give the proofs of Theorem \ref{thm kaehler} and Theorem \ref{thm n-1 flow}.

\subsection{Theorem \ref{thm kaehler}}
\begin{proof}
The first observation is that, when $k=1$, the inequalities in the numerical conditions are just the right hand side in weak transcendental holomorphic Morse inequalities. Recall that Demailly's conjecture on weak transcendental holomorphic Morse inequalities (see e.g. \cite[Conjecture 10.1]{BDPP13}) is stated as following:\\

\noindent
\emph{Let $X$ be a compact complex manifold of dimension $n$, and let $\gamma, \beta$ be two nef classes over $X$. Then we have $$\vol{(\gamma-\beta)}\geq \gamma^n -n \gamma^{n-1}\cdot \beta.$$
In particular, $\gamma^n -n \gamma^{n-1}\cdot \beta >0$ implies the class $\gamma-\beta$ is big, that is, $\gamma-\beta$ contains a K\"ahler current.}
\\

\noindent
Note that the last statement has been proved for K\"ahler manifolds \cite{Pop14} (see also \cite{Xia13}), that is, if $X$ is a compact K\"ahler manifold then $\gamma^n -n \gamma^{n-1}\cdot \beta >0$ implies there exists a K\"ahler current in the class $\gamma-\beta$.

We apply this bigness criterion to the classes $\{\omega\}$ and $ \{\alpha\}$,
then the numerical condition (\ref{eq numerical positive}) implies $\{\omega-\alpha\}_{|_V}$ is a big class on every proper irreducible subvariety $V$. More precisely, if $V$ is singular then by some resolution of singularities we have a proper modification $\pi : \widehat{V}\rightarrow V$ with $\widehat{V}$ smooth, and by (\ref{eq numerical positive}) we know
\begin{align*}
\pi^*\{\omega\}_{|_V} ^p- p \pi^*\{\omega\}_{|_V} ^{p-1}\cdot \pi^*\{\alpha\}_{|_V}>0,
\end{align*}
thus the class $\pi^*\{\omega-\alpha\}_{|_V}$ contains a K\"ahler current over $\widehat{V}$. So by the push-forward map $\pi_*$ we obtain that the class $\{\omega-\alpha\}_{|_V}$ is big over $V$.

In particular, by (\ref{eq over X}) and (\ref{eq numerical positive}) the restriction of the class $\{\omega-(1-\epsilon)\alpha\}$ is big on every irreducible subvariety (including $X$ itself) for any sufficiently small $\epsilon>0$.

We claim this yields $\{\omega-(1-\epsilon)\alpha\}$ is a K\"ahler class over $X$ for any $\epsilon>0$ small.
Indeed, our proof implies the following fact.
\begin{itemize}
  \item Assume $\beta$ is a big class over a compact complex manifold (or compact complex space) and its restriction to every irreducible subvariety is also big, then $\beta$ is a K\"ahler class over $X$.
\end{itemize}
To this end, we will argue by induction on the dimension of $X$. If $X$ is a compact complex curve, then this is obvious. For the general case, we need a result of Mihai P{\u{a}}un (see \cite{paun1998effectivite, paun98thesis}):\\

\noindent
\emph{Let $X$ be a compact complex manifold (or compact complex space), and let $\beta=\{T\}$ be the cohomology class of a K\"ahler current $T$ over $X$. Then $\beta$ is a K\"ahler class over $X$ if and only if the restriction $\beta_{|_Z}$ is a K\"ahler class on every irreducible component $Z$ of the Lelong sublevel set
$E_c (T)$.}\\

\noindent
As $\{\omega-(1-\epsilon)\alpha\}$ is a big class on $X$, by Demailly's regularization theorem \cite{Dem92} we can choose a K\"ahler current
$T\in \{\omega-(1-\epsilon)\alpha\}$ such that $T$ has analytic singularities on $X$.
Then the singularities of $T$ are just the Lelong sublevel set $E_c(T)$ for some positive constant $c$. For every irreducible component $Z$ of $E_c (T)$, by (\ref{eq numerical positive}) the restriction $\{\omega-(1-\epsilon)\alpha\}_{|_Z}$ is a big class. After resolution of singularities of $Z$ if necessary, we obtain a K\"ahler current $T_Z \in \{\omega-(1-\epsilon)\alpha\}_{|_Z}$ over $Z$ with its analytic singularities contained in a proper subvariety of $Z$, and for every irreducible subvariety $V\subseteq Z$ the restriction $\{\omega-(1-\epsilon)\alpha\}_{|_V}$ is also a big class. By induction on the dimension, we get that $\{\omega-(1-\epsilon)\alpha\}_{|_Z}$ is a K\"ahler class over $Z$. So the above result of \cite{paun1998effectivite, paun98thesis} implies $\{\omega-(1-\epsilon)\alpha\}$ is a K\"ahler class over $X$, finishing the proof our claim.

By the arbitrariness of $\epsilon>0$, we get $\{\omega-\alpha\}$ is a nef class on $X$.
Next we prove $\{\omega-\alpha\}$ is a big class. By \cite[Theorem 2.12]{DP04}, we only need to show
$$
\vol(\{\omega-\alpha\})=\int_X (\omega-\alpha)^n>0.
$$
Since $\{\omega-\alpha\}$ is nef, we can compute the derivative of the function $\vol(\omega-t\alpha)$ for any $t\in [0,1)$. Thus we have
\begin{align*}
\vol(\{\omega\}-\{\alpha\})-\vol(\{\omega\})&=\int_0 ^1 \frac{d}{dt}\vol(\{\omega\}-t\{\alpha\})dt\\
&=-\int_0 ^1 n\{\omega-t\alpha\}^{n-1}\cdot \{\alpha\} dt,
\end{align*}
which implies
\begin{align*}
\vol(\{\omega\}-\{\alpha\})&=\vol(\{\omega\})-\int_0 ^1 n\{\omega-t\alpha\}^{n-1}\cdot \{\alpha\} dt\\
&\geq \int_0 ^1 n(\{\omega\}^{n-1}-\{\omega-t\alpha\}^{n-1})\cdot \{\alpha\} dt.
\end{align*}
Here the last line follows from the equality (\ref{eq over X}). Since $\omega, \alpha$ are K\"ahler metrics, this shows $\vol(\{\omega-\alpha\})>0$. Thus $\{\omega-\alpha\}$ is a big and nef class on $X$ with its restriction to every irreducible subvariety being big and nef. By the arguments before, we know $\{\omega-\alpha\}$ must be a K\"ahler class.

Finally, we give an alternative proof of the fact that the class $\{\omega-\alpha\}$ is nef using the main result of \cite{collins2013null-locus} instead of using \cite{paun1998effectivite, paun98thesis}. (I would like to thank Tristan C. Collins who pointed out this to me.)
Since $\{\omega\}$ is a K\"ahler class, the class $\{\omega-t\alpha\}$ is K\"ahler for $t>0$ small. Let $s$ be the largest number such that $\{\omega-s\alpha\}$ is nef. We prove that $s\geq 1$. Otherwise, suppose $s<1$. Then by the numerical conditions (\ref{eq over X}) and (\ref{eq numerical positive}), the bigness criterion given by transcendental holomorphic Morse inequalities implies that the class $\{\omega-s\alpha\}$ is big if $s<1$, and furthermore, this holds for all irreducible subvarieties in $X$.  Thus $\{\omega-s\alpha\}$ is big and nef on every irreducible subvariety $V$ in $X$.
This means the null locus of the big and nef class $\{\omega-s\alpha\}$ is empty, and then the main result of \cite{collins2013null-locus} implies that $\{\omega-s\alpha\}$ is a K\"ahler class. This contradicts with the definition of $s$, so we get $s\geq 1$, or equivalently, $\{\omega-\alpha\}$ must be a nef class. Then by the estimate of the volume $\vol(\{\omega-\alpha\})$ as above, we know $\{\omega-\alpha\}$ is also big and nef over every irreducible subvariety of $X$. By applying \cite{collins2013null-locus} again, this proves that $\{\omega-\alpha\}$ must be a K\"ahler class.
\end{proof}

\begin{remark}
If $X$ is a smooth projective variety of dimension $n$ and $\{\omega\}$ and $\{\alpha\}$ are the first Chern classes of holomorphic line bundles, then the nefness of the class $\{\omega-\alpha\}$ just follows from Kleiman's ampleness criterion, since the numerical condition (\ref{eq numerical positive}) for $p=1$ implies the divisor class $\{\omega-\alpha\}$ has non-negative intersection against every irreducible curve.
\end{remark}

\subsection{Theorem \ref{thm n-1 flow}}
\label{section thm n-1 flow}
Next we give the proof of Theorem \ref{thm n-1 flow}.

\begin{proof}
The proof mainly depends on Boucksom's divisorial Zariski decomposition for pseudoeffective $(1,1)$-classes \cite{Bou04} and the bigness criterion for the difference of two movable $(n-1, n-1)$-classes \cite{xiao2014movable}.

Through a sufficiently small perturbation of the K\"ahler metric $\alpha$, e.g. replace $\alpha$ by
$$\alpha_\epsilon = (1-\epsilon)\alpha$$
with $\epsilon\in (0, 1)$, we can obtain that the inequality in (\ref{eq over X}) is strict for the classes $\{\omega\}$ and $\{\alpha_\epsilon\}$. We claim that in this case the $(n-1, n-1)$-class $\{\omega^{n-1}- \alpha_\epsilon ^{n-1}\}$ has nonnegative intersections with all pseudoeffective $(1,1)$-classes. Then let $\epsilon$ tends to zero, we conclude the desired result for the class
$\{\omega^{n-1}- \alpha ^{n-1}\}$. Thus we can assume the inequality in (\ref{eq over X})is strict for the classes $\{\omega\}$ and $\{\alpha\}$ at the beginning.

Let $\beta$ be a pseudoeffective $(1,1)$-class over $X$. 
By \cite[Section 3]{Bou04}, $\beta$ admits a divisorial Zariski decomposition
\begin{align*}
\beta= Z(\beta)+N(\beta).
\end{align*}
Note that $N(\beta)$ is the class of some effective divisor (may be zero) and $Z(\beta)$ is a modified nef class. In particular, we have
\begin{align}\label{eq positive negativepart}
\{\omega^{n-1}-\alpha^{n-1}\}\cdot N(\beta)\geq 0.
\end{align}
For any $\delta>0$, we have
\begin{align*}
 Z(\beta)+\delta \{\omega\} = \pi_* \{\widehat{\omega}\}
\end{align*}
for some modification $\pi: \widehat{X}\rightarrow X$ and some K\"ahler metric $\widehat{\omega}$ on $\widehat{X}$ (see \cite[Proposition 2.3]{Bou04}).

By our assumption on (\ref{eq over X}), we have
\begin{align}\label{eq pullback n-1 positive}
 \int_{\widehat{X}} \pi^*\omega^n -n \pi^*\omega \wedge \pi^*\alpha^{n-1}>0.
\end{align}
By \cite[Theorem 3.3]{xiao2014movable} (or \cite[Remark 3.1]{Xia13}), the inequality (\ref{eq pullback n-1 positive}) implies that the class $\{\pi^* \omega^{n-1}-\pi^* \alpha^{n-1}\}$ contains a strictly positive $(n-1, n-1)$-current. This implies
\begin{align*}
&\{\omega^{n-1}-\alpha^{n-1}\}\cdot (Z(\beta)+\delta \{\omega\})\\
&=\{\omega^{n-1}-\alpha^{n-1}\}\cdot \pi_*\{\widehat{\omega}\}\\
&=\pi^*\{\omega^{n-1}-\alpha^{n-1}\}\cdot \{\widehat{\omega}\}\\
&>0.
\end{align*}
By the arbitrariness of $\delta$, we get $\{\omega^{n-1}-\alpha^{n-1}\}\cdot Z(\beta) \geq 0$. With (\ref{eq positive negativepart}), we show that
$$
\{\omega^{n-1}-\alpha^{n-1}\}\cdot \beta \geq 0.
$$
Since $\beta$ can be any pseudoeffective $(1,1)$-class, this implies $\{\omega^{n-1}-\alpha^{n-1}\}$ lies in the closure of the Gauduchon cone by \cite[Proposition 2.1]{xiao2015characterizing} (see also \cite[Lemma 3.3]{lamari99lemma}).
\end{proof}

\begin{remark}
We expect $\{\omega^{n-1}-\alpha^{n-1}\}$ should have strictly positive intersection numbers with nonzero pseudoeffective $(1,1)$-classes. To show this, one only need to verify this for modified nef classes.
\end{remark}

\begin{remark}
Let $X$ be a smooth projective variety, and assume $\{\omega^{n-1}- \alpha^{n-1}\}$ is a curve class. Then the numerical condition (\ref{eq numerical positive}) in Theorem \ref{thm n-1 flow} implies that $\{\omega^{n-1}- \alpha^{n-1}\}$ is a movable class by \cite[Theorem 2.2]{BDPP13}.
\end{remark}

\section{Further discussions}
\label{section discussions}

In analogue with Theorem \ref{thm kaehler} and Theorem \ref{thm n-1 flow}, one would like to prove similar positivity of the class $\{\omega^k - \alpha^k\}$. To generalize our results in this direction, one can apply \cite[Remark 3.1]{Xia13}. By \cite[Remark 3.1]{Xia13}, we know that the condition
\begin{align*}
\label{eq numerical positive}
\int_V \omega^p - \frac{p!}{k!(p-k)!} \omega^{p-k}\wedge \alpha^k >0
\end{align*}
implies that the class $\{\omega^k - \alpha^k\}_{|V}$ contains a strictly positive $(k,k)$-current
over every irreducible subvariety $V$ of dimension $p$ with
$k< p \leq n-1$.
However, the difficulties appear as we know little about the singularities of positive $(k,k)$-currents for $k>1$. We have no analogues of Demailly's regularization theorem for such currents.

Inspired by the prediction of Conjecture \ref{conj LS15 sigma k flow}, we propose the following question on the positivity of positive $(k, k)$-currents.

\begin{question}
\label{question kk current}
Let $X$ be a compact K\"ahler manifold (or general compact complex manifold) of dimension $n$. Let $\Omega\in H^{k,k}(X, \mathbb{R})$ be a big $(k,k)$-class, i.e. it can be represented by a strictly positive $(k,k)$-current over $X$. Assume the restriction class $\Omega_{|V}$ is also big over every irreducible subvariety $V$ with $k\leq \dim V \leq n-1$, then does $\Omega$ contain a smooth strictly positive $(k,k)$-form in its Bott-Chern class? Or does $\Omega$ contain a strictly positive $(k,k)$-current with analytic singularities of codimension at least $n-k+1$ in its Bott-Chern class?
\end{question}

\bibliography{reference}

\providecommand{\bysame}{\leavevmode\hbox to3em{\hrulefill}\thinspace}
\providecommand{\MR}{\relax\ifhmode\unskip\space\fi MR }
\providecommand{\MRhref}[2]{%
  \href{http://www.ams.org/mathscinet-getitem?mr=#1}{#2}
}
\providecommand{\href}[2]{#2}
\begin{thebibliography}{BDPP13}

\bibitem[BDPP13]{BDPP13}
S{\'e}bastien Boucksom, Jean-Pierre Demailly, Mihai P{\u{a}}un, and Thomas
  Peternell, \emph{The pseudo-effective cone of a compact {K}\"ahler manifold
  and varieties of negative {K}odaira dimension}, J. Algebraic Geom.
  \textbf{22} (2013), no.~2, 201--248. \MR{3019449}

\bibitem[Bou04]{Bou04}
S{\'e}bastien Boucksom, \emph{Divisorial {Z}ariski decompositions on compact
  complex manifolds}, Ann. Sci. \'Ecole Norm. Sup. (4) \textbf{37} (2004),
  no.~1, 45--76. \MR{2050205 (2005i:32018)}

\bibitem[Che00]{Chen00mabuchi_Jflow}
Xiuxiong Chen, \emph{On the lower bound of the {M}abuchi energy and its
  application}, Internat. Math. Res. Notices (2000), no.~12, 607--623.
  \MR{1772078 (2001f:32042)}

\bibitem[Che04]{Chen04Jflow}
\bysame, \emph{A new parabolic flow in {K}\"ahler manifolds}, Comm. Anal. Geom.
  \textbf{12} (2004), no.~4, 837--852. \MR{2104078 (2005h:53116)}

\bibitem[CS14]{collins2014convergence-Jflow}
Tristan~C. Collins and G{\'a}bor Sz{\'e}kelyhidi, \emph{Convergence of the
  $j$-flow on toric manifolds}, arXiv preprint arXiv:1412.4809 (2014).

\bibitem[CT13]{collins2013null-locus}
Tristan~C. Collins and Valentino Tosatti, \emph{K\"ahler currents and null
  loci}, arXiv preprint arXiv:1304.5216 (2013).

\bibitem[Dem92]{Dem92}
Jean-Pierre Demailly, \emph{Regularization of closed positive currents and
  intersection theory}, J. Algebraic Geom. \textbf{1} (1992), no.~3, 361--409.
  \MR{1158622 (93e:32015)}

\bibitem[Don99]{Donaldson99Jflow_momentmap}
Simon~K. Donaldson, \emph{Moment maps and diffeomorphisms}, Asian J. Math.
  \textbf{3} (1999), no.~1, 1--15, Sir Michael Atiyah: a great mathematician of
  the twentieth century. \MR{1701920 (2001a:53122)}

\bibitem[DP04]{DP04}
Jean-Pierre Demailly and Mihai P{\u{a}}un, \emph{Numerical characterization of
  the {K}\"ahler cone of a compact {K}\"ahler manifold}, Ann. of Math. (2)
  \textbf{159} (2004), no.~3, 1247--1274. \MR{2113021 (2005i:32020)}

\bibitem[FLM11]{maxinan11sigmaflow}
Hao Fang, Mijia Lai, and Xinan Ma, \emph{On a class of fully nonlinear flows in
  {K}\"ahler geometry}, J. Reine Angew. Math. \textbf{653} (2011), 189--220.
  \MR{2794631 (2012g:53134)}

\bibitem[Lam99]{lamari99lemma}
Ahc{\`e}ne Lamari, \emph{Courants k\"ahl\'eriens et surfaces compactes}, Ann.
  Inst. Fourier (Grenoble) \textbf{49} (1999), no.~1, vii, x, 263--285.
  \MR{1688140 (2000d:32034)}

\bibitem[LS15]{lejmisze15Jflow}
Mehdi Lejmi and G{\'a}bor Sz{\'e}kelyhidi, \emph{The {J}-flow and stability},
  Adv. Math. \textbf{274} (2015), 404--431. \MR{3318155}

\bibitem[P{\u{a}}u98a]{paun98thesis}
Mihai P{\u{a}}un, \emph{Fibr\'{e} en droites num\'{e}riquement effectifs et
  vari\'{e}t\'{e}s k\"ahl\'{e}riennes compactes \`{a} courbure de {R}icci nef},
  Ph.D. thesis, Universit{\'e} Joseph-Fourier-Grenoble I, 1998.

\bibitem[P{\u{a}}u98b]{paun1998effectivite}
\bysame, \emph{Sur l'effectivit{\'e} num{\'e}rique des images inverses de
  fibr{\'e}s en droites}, Mathematische Annalen \textbf{310} (1998), no.~3,
  411--421.

\bibitem[Pop14]{Pop14}
Dan Popovici, \emph{An observation relative to a paper by {J}. {X}iao}, arXiv
  preprint arXiv:1405.2518 (2014).

\bibitem[SW08]{songweinkove08Jflow}
Jian Song and Ben Weinkove, \emph{On the convergence and singularities of the
  {$J$}-flow with applications to the {M}abuchi energy}, Comm. Pure Appl. Math.
  \textbf{61} (2008), no.~2, 210--229. \MR{2368374 (2009a:32038)}

\bibitem[Xia13]{Xia13}
Jian Xiao, \emph{Weak transcendental holomorphic {M}orse inequalities on
  compact {K}\"ahler manifolds}, arXiv preprint arXiv:1308.2878, to appear in
  Ann. Inst. Fourier (Grenoble) (2013).

\bibitem[Xia14]{xiao2014movable}
\bysame, \emph{Movable intersection and bigness criterion}, arXiv preprint
  arXiv:1405.1582 (2014).

\bibitem[Xia15]{xiao2015characterizing}
\bysame, \emph{Characterizing volume via cone duality}, arXiv preprint
  arXiv:1502.06450 (2015).

\end{thebibliography}
\bibliographystyle{amsalpha}

\noindent
\textsc{Institute of Mathematics, Fudan University, 200433 Shanghai, China}\\

\noindent
\textsc{Current address:}\\
\textsc{Institut Fourier, Universit\'{e} Joseph Fourier, 38402 Saint-Martin d'H\`{e}res, France}
\verb"Email: jian.xiao@ujf-grenoble.fr"

\end{document}